\documentstyle{article}
\begin{document}
\renewcommand{\baselinestretch}{1.2}
\begin{center}
{\Large\bf {Jordan higher all-derivable points in triangular
algebras }} {\footnote{This work is supported by the National
Natural Science Foundation of China (No 10771191) }}\vspace{0.4cm}
 \\
{\bf Jinping Zhao \footnote{E-mail address: zjphyx@sohu.com}, Jun zhu{\footnote{E-mail address: zhu$\_$gjun@yahoo.com.cn}}}\\
\vspace{0.3cm} \small{Institute of Mathematics, Hangzhou Dianzi
University, Hangzhou 310018, People's Republic of China}
\end{center}
\underline{~~~~~~~~~~~~~~~~~~~~~~~~~~~~~~~~~~~~~~~~~~~~~~~~~~~~~~~~~~~~~~~~~
~~~~~~~~~~~~~~~~~~~~~~~~~~~~~~~~~~~~~~~~~~~~~~~~~~~~~~~~~~~~~~~}
\vspace{0.1cm}\ \vspace{0.2cm} {\bf{Abstract}}^^L {Let
${\mathcal{T}}$ be a triangular algebra. We say that $D=\{D_{n}:
n\in N\}\subseteq L({\mathcal{T}})$ is a Jordan higher derivable
mapping at $G$ if
$D_{n}(ST+TS)=\sum_{i+j=n}(D_{i}(S)D_{j}(T)+D_{i}(T)D_{j}(S))$ for
any $S,T\in {\mathcal{T}}$ with $ST=G$. An element $G\in
{\mathcal{T}}$ is called a Jordan higher all-derivable point of
${\mathcal{T}}$ if every Jordan higher derivable linear mapping
$D=\{D_{n}\}_{n\in N}$ at $G$ is a higher derivation. In this paper,
under some mild conditions on ${\mathcal{T}}$, we prove that some
elements of ${\mathcal{T}}$ are Jordan higher all-derivable points.
This extends some results in [6] to the case of Jordan higher
derivations. \vspace{0.2cm}
\\
\vspace{0.2cm} {\it{AMS Classification}}:16W25, 47B47}\\
\vspace{0.2cm}
\underline{~~~~~~~~~~~~~~~~~~~~~~~~~~~~~~~~~~~~~~~~~~~~~~~~~~~~~~~~~~~~~~~~~
~~~~~~~~~~~~~~~~~~~~~~~~~~~~~~~~~~~~~~~~~~~~~~~~~~~~~~~~~~~~~}
\vspace{0.1cm}\

\section*{1. Introduction and preliminaries}
~

Let $\mathcal{A}$  be a ring (or algebra) with the unit $I$. An
additive linear  mapping $\delta$ from $\mathcal{A}$ into itself is
called a derivation if $\delta(ST)=\delta(S)T+S\delta(T)$ for any
$S, T\in \mathcal{A}$ and is said to be a Jordan derivation if
$\delta(ST+TS)=\delta(S)T+S\delta(T)+\delta(T)S+T\delta(S)$ for any
$S, T\in \mathcal{A}$. We say that a mapping $\delta$ is Jordan
derivable at a given point $G\in \mathcal{A}$ if
$\delta(ST+TS)=\delta(S)T+S\delta(T)+\delta(T)S+T\delta(S)$ for any
$S, T\in \mathcal{A}$ with $ST=G$, and $G$ is called a Jordan
all-derivable point of $\mathcal{A}$ if every Jordan derivable
mapping at $G$ is a derivation. We say that $D=\{D_{n}\}\subseteq
L({\mathcal{A}})$ is a Jordan higher derivable mapping at $G$ if
$D_{n}(ST+TS)=\sum_{i+j=n}(D_{i}(S)D_{j}(T)+D_{i}(T)D_{j}(S))$ for
any $S,T\in {\mathcal{A}}$ with $ST=G$. An element $G\in
{\mathcal{A}}$ is called a Jordan higher all-derivable point of
${\mathcal{A}}$ if every Jordan higher derivable linear mapping
$D=\{D_{n}\}$ at $G$ is a higher derivation. There have been a
number of papers on the study of conditions under which derivations
of operator algebras can be completely determined by the action on
some sets of operators. In [3], W. Jing showed that I is a Jordan
all-derivable point of $\mathcal{B}(\mathcal{H})$ with $\mathcal{H}$
is a Hilbert space. In [7], J. Zhu proved that every invertible
operator in nest algebra is an all-derivable point in the strong
operator topology. Also it was showed that every element in
 the algebra of all upper triangular matrices is a Jordan
all-derivable point by Z. Sha and J. Zhu in [6].

With the development of derivation, higher derivation has attracted
much attention of  mathematicians as an active subject of research
in algebras. In [4] Z. Xiao and F. Wei showed that any Jordan higher
derivation on a triangular algebra is a higher derivation. In this
paper we will extend the conclusion of [6] to the case of Jordan
higher derivations.

Let $\mathcal{A}$ and $\mathcal{B}$ be two unital rings (or
algebras) with the unit $I_{1}$, $I_{2}$, and  $\mathcal{M}$ be a
 unital ($\mathcal{A}$, $\mathcal{B}$)-bimodule, which is faithfull as a left  $\mathcal{A}$-module and as a right
 $\mathcal{B}$-module. The ring(or algebra)

$${\mathcal{T}}=\{{\left [\begin{array}{ccc}
  a & m \\
  0 & b \\
\end{array}\right ]}: a\in
{\mathcal{A}}, m\in {\mathcal{M}}, b\in {\mathcal{B}}\},$$ under the
usual matrix operations is said to be a triangular algebra. We
mainly proved that $0$ and $\left [\begin{array}{ccc}
  I_{1} & X_{0} \\
  0 & I_{2} \\
\end{array}\right ]$ are Jordan higher  all-derivable points for
any given point $X_{0}\in \mathcal{M}$.

\section*{2. Jordan higher all-derivable points in ring algebras}
~

In this section, we always assume that the characteristics of
$\mathcal{A}$ and $\mathcal{B}$ are not $2$ and $3$, and for any
$X\in \mathcal{A}$, $Y\in \mathcal{B}$, there are some integers
$n_{1}$, $n_{2}$ such that $n_{1}I_{1}-X$ and $n_{2}I_{2}-Y$ are
invertible. The following two theorems are the main results in this paper.\\
~\\ {\bf{Theorem 2.1}} {\it Let $D=(D_{n})_{n\in N}$ be a family of
additive linear mappings on $\mathcal{T}$ that
${D}_{0}=iD_{\mathcal{T}}$ (identical mapping on ${\mathcal{T}}$).
If $D$ is Jordan
higher derivable at $0$, then $D$ is a higher derivation.}\\
{\bf Proof.} For any $T=\left [\begin{array}{cc}
  X & Y \\
  0 & Z \\
\end{array}\right ]\in {\mathcal{T}}$, we can write $$D_{n}(\left [\begin{array}{cc}
  X & Y \\
  0 & Z \\
\end{array}\right ])=\left [\begin{array}{cc}
  \delta_{n}^{11}(X)+\varphi_{n}^{11}(Y)+\tau_{n}^{11}(Z) & \delta_{n}^{12}(X)+\varphi_{n}^{12}(Y)+\tau_{n}^{12}(Z)  \\
  0 & \delta_{n}^{22}(X)+\varphi_{n}^{22}(Y)+\tau_{n}^{22}(Z)  \\
\end{array}\right ],$$ where $\delta_{n}^{ij}:{\mathcal{A}}\rightarrow {\mathcal{A}}_{ij}$,
$\varphi_{n}^{ij}:{\mathcal{M}}\rightarrow {\mathcal{A}}_{{ij}}$,
$\tau_{n}^{ij}:{\mathcal{B}}\rightarrow {\mathcal{A}}_{ij}$, $1\leq
i\leq j\leq 2$ are additive maps with
${\mathcal{A}}_{11}=\mathcal{A}$, ${\mathcal{A}}_{12}=\mathcal{M}$,
${\mathcal{A}}_{22}=\mathcal{B}$. It follows from the fact
${D}_{0}=iD_{\mathcal{T}}$ that when $i=j=1$,
$\delta_{0}^{ij}=i\delta_{\mathcal{A}}$, else $\delta_{0}^{ij}=0$;
when $i=1,j=2$, $\varphi_{0}^{ij}=i\varphi_{\mathcal{M}}$, else
$\varphi_{0}^{ij}=0$;
 when $i=j=2$
, $\tau_{0}^{ij}=i\tau_{\mathcal{B}}$, else $\tau_{0}^{ij}=0$.

We set $S=\left [\begin{array}{cc}
  0 & W \\
  0 & 0 \\
\end{array}\right ]$ and $T=\left
[\begin{array}{cc}
  X & 0 \\
  0 & 0 \\
\end{array}\right ]$ for every $X\in {\mathcal{A}}$, $W\in {\mathcal{M}}$. Then $ST=0$ and $TS=\left
[\begin{array}{cc}
  0 & XW \\
  0 & 0 \\
\end{array}\right ]$. So $$\begin{array}{rcl}& &\left
[\begin{array}{ccc}
  \varphi_{n}^{11}(XW) & \varphi_{n}^{12}(XW) \\
  0 & \varphi_{n}^{22}(XW) \\
\end{array}\right
]=D_{n}(ST+TS)=\sum\limits_{i+j=n}(D_{i}(S)D_{j}(T)+D_{i}(T)D_{j}(S))
\\\\&=& \sum\limits_{i+j=n}(\left [\begin{array}{cc}
 \varphi_{i}^{11}(W) & \varphi_{i}^{12}(W) \\
  0 & \varphi_{i}^{22}(W) \\
\end{array}\right ]\left [\begin{array}{ccc}
  \delta_{j}^{11}(X) & \delta_{j}^{12}(X) \\
  0 & \delta_{j}^{22}(X) \\
\end{array}\right ]\\&&+\left [\begin{array}{ccc}
  \delta_{i}^{11}(X) & \delta_{i}^{12}(X) \\
  0 & \delta_{i}^{22}(X) \\
\end{array}\right ]\left [\begin{array}{cc}
 \varphi_{j}^{11}(W) & \varphi_{j}^{12}(W) \\
  0 & \varphi_{j}^{22}(W) \\
\end{array}\right ])
\\\\&=& \sum\limits_{i+j=n}\left [\begin{array}{cc}
 \varphi_{i}^{11}(W)\delta_{j}^{11}(X)+\delta_{i}^{11}(X)\varphi_{j}^{11}(W) &
 \varphi_{i}^{11}(W)\delta_{j}^{12}(X)+\delta_{i}^{11}(X)\varphi_{j}^{12}(W)
  \\~&+\varphi_{i}^{12}(W)\delta_{j}^{22}(X)+\delta_{i}^{12}(X)\varphi_{j}^{22}(W)\\~&~\\
  0 & \varphi_{i}^{22}(W)\delta_{j}^{22}(X)+\delta_{i}^{22}(X)\varphi_{j}^{22}(W) \\
\end{array}\right ].\\
\end{array}$$ This implies that
\begin{equation}
 \varphi_{n}^{11}(XW)=\sum\limits_{i+j=n}
 (\varphi_{i}^{11}(W)\delta_{j}^{11}(X)+\delta_{i}^{11}(X)\varphi_{j}^{11}(W)),
\end{equation}
\begin{equation}\
 \varphi_{n}^{12}(XW)=\sum\limits_{i+j=n}(\varphi_{i}^{11}(W)\delta_{j}^{12}(X)+\delta_{i}^{11}(X)\varphi_{j}^{12}(W)
  +\varphi_{i}^{12}(W)\delta_{j}^{22}(X)+\delta_{i}^{12}(X)\varphi_{j}^{22}(W)),
\end{equation}
and
\begin{equation}
\varphi_{n}^{22}(XW)=\sum\limits_{i+j=n}(\varphi_{i}^{22}(W)\delta_{j}^{22}(X)+\delta_{i}^{22}(X)\varphi_{j}^{22}(W))
\end{equation}
for any $X\in {\mathcal{A}}$, $W\in {\mathcal{M}}$.  One obtains
that
\begin{equation}
 \varphi_{n}^{11}(W)=\sum\limits_{i+j=n}
 (\varphi_{i}^{11}(W)\delta_{j}^{11}(I_{1})+\delta_{i}^{11}(I_{1})\varphi_{j}^{11}(W)),
\end{equation}
\begin{equation}
\varphi_{n}^{22}(W)=\sum\limits_{i+j=n}(\varphi_{i}^{22}(W)\delta_{j}^{22}(I_{1})+\delta_{i}^{22}(I_{1})\varphi_{j}^{22}(W))
\end{equation}
by taking $X=I_{1}$ in Eq. (1) and Eq. (3). Now we prove the fact
that $\varphi_{n}^{11}(W)=0$ and $\varphi_{n}^{22}(W)=0$ by
induction on n. When $n=0$, it is easily verified that
$\varphi_{0}^{11}(W)=0$ and $\varphi_{0}^{22}(W)=0$ from the
characterizations of $\varphi_{0}^{11}$ and $\varphi_{0}^{22}$. When
$n=1$, $\varphi_{1}^{11}(W)=0$ and $\varphi_{1}^{22}(W)=0$ can be
obtained by the proof in [6, Theorem 2.1]. We assume that
$\varphi_{m}^{11}(W)=0$ and $\varphi_{m}^{22}(W)=0$ for all $1\leq
m<n$. In fact, by the Eq. (4) and
$\delta_{0}^{11}=i\delta_{\mathcal{A}}$, we have
$\varphi_{n}^{11}(W)=\varphi_{n}^{11}(W)+\varphi_{n}^{11}(W)=2\varphi_{n}^{11}(W)$.
Thus $\varphi_{n}^{11}(W)=0$. Similarly combining Eq. (5) with the
fact that $\delta_{0}^{22}=0$, we can get $\varphi_{n}^{22}(W)=0$
for any $W\in M$ and $n\in N$.
 For any $X\in
\mathcal{A}$, $W\in \mathcal{M}$ and $Y\in \mathcal{B}$, setting
$S=\left [\begin{array}{cc}
  0 & W \\
  0 & Y \\
\end{array}\right ]$ and $T=\left
[\begin{array}{cc}
  X & 0 \\
  0 & 0 \\
\end{array}\right ]$, then $ST=0$, $TS=\left [\begin{array}{cc}
  0 & XW \\
  0 & 0 \\
\end{array}\right ]$. One gets $$\begin{array}{rcl}& &\left
[\begin{array}{ccc}
  0 & \varphi_{n}^{12}(XW) \\
  0 & 0 \\
\end{array}\right
]=D_{n}(ST+TS)=\sum\limits_{i+j=n}(D_{i}(S)D_{j}(T)+D_{i}(T)D_{j}(S))
\\\\&=&\sum\limits_{i+j=n}(\left [\begin{array}{cc}
\tau_{i}^{11}(Y) & \varphi_{i}^{12}(W)
+\tau_{i}^{12}(Y) \\
  0 & \tau_{i}^{22}(Y) \\
\end{array}\right ]\left [\begin{array}{ccc}
  \delta_{j}^{11}(X) & \delta_{j}^{12}(X) \\
  0 & \delta_{j}^{22}(X) \\
\end{array}\right ]\\\\&&+\left [\begin{array}{ccc}
  \delta_{i}^{11}(X) & \delta_{i}^{12}(X) \\
  0 & \delta_{i}^{22}(X) \\
\end{array}\right ]\left [\begin{array}{cc}
 \tau_{j}^{11}(Y) & \varphi_{j}^{12}(W)+\tau_{j}^{12}(Y) \\
  0 & \tau_{j}^{22}(Y) \\
\end{array}\right ]).
\end{array}$$
Hence the following three equations hold
\begin{equation}
 \sum\limits_{i+j=n}
 (\tau_{i}^{11}(Y)\delta_{j}^{11}(X)+\delta_{i}^{11}(X)\tau_{j}^{11}(Y))=0,
\end{equation}
\begin{equation}
\sum\limits_{i+j=n}(\tau_{i}^{22}(Y)\delta_{j}^{22}(X)+\delta_{i}^{22}(X)\tau_{j}^{22}(Y))=0,
\end{equation}
\begin{equation}\begin{array}{rcl}
\varphi_{n}^{12}(XW)&=&
\sum\limits_{i+j=n}(\tau_{i}^{11}(Y)\delta_{j}^{12}(X)+\varphi_{i}^{12}(W)\delta_{j}^{22}(X)+\tau_{i}^{12}(Y)\delta_{j}^{22}(X)
\\\\&&+\delta_{i}^{11}(X)\varphi_{j}^{12}(W)+\delta_{i}^{11}(X)\tau_{j}^{12}(Y)+\delta_{i}^{12}(X)\tau_{j}^{22}(Y)
\end{array}
\end{equation}
for any $X \in {\mathcal{A}}$, $W \in {\mathcal{M}}$. One can see
that
\begin{equation}
 \sum\limits_{i+j=n}
 (\tau_{i}^{11}(Y)\delta_{j}^{11}(I_{1})+\delta_{i}^{11}(I_{1})\tau_{j}^{11}(Y))=0
\end{equation} by taking $X=I_{1}$ in Eq. (6). Using Eq. (9) and
induction, one has $\tau_{n}^{11}(Y)=0$ for every $n\in N$.
Similarly taking $Y=I_{2}$ in Eq. (7), by inducting and using the
fact that $\tau_{0}^{22}(Y)=0$, we get $\delta_{n}^{22}(X)=0$ for
every $n\in N$ and $X\in \mathcal{A}$.

We can obtain that
\begin{equation}
\sum\limits_{i+j=n}(\delta_{i}^{11}(X)\tau_{j}^{12}(Y)+\delta_{i}^{12}(X)\tau_{j}^{22}(Y))=0
\end{equation}
\\
by $\delta^{22}_{i}(X)=0$, $\tau^{11}_{i}(Y)=0$ and taking $W=0$ in
Eq. (8).

By Eq. (2) and the fact that $\delta_{n}^{22}(X)=0$,
$\varphi_{n}^{11}(W)=0$, $\varphi_{n}^{22}(W)=0$ and
$\varphi_{0}^{12}=i\varphi_{\mathcal{M}}$, we have
\begin{equation}
\varphi_{n}^{12}(XW)= \sum\limits_{i+j=n}
\delta_{i}^{11}(X)\varphi_{j}^{12}(W).
\end{equation}

We claim that $\delta=\{\delta^{11}_{n}: n\in N\}$ is a higher
derivation on $\mathcal{A}$. In fact, we know that $\delta_{1}$ is a
derivation by Theorem 2.1 in [6]. It follows that
$\delta_{1}^{11}(X_{1}X_{2})=\delta_{1}^{11}(X_{1})X_{2}+X_{1}\delta_{1}^{11}(X_{2})$
for any $X_{1}, X_{2}$ in $\mathcal{A}$. Now we assume that
$\delta_{m}^{11}(X_{1}X_{2})= \sum\limits_{i+j=m}
\delta_{i}^{11}(X_{1})\delta_{j}^{11}(X_{2})$ for any $1\leq m<n$
with $m\in N$. Summing up Eq. (11) and
$\varphi_{0}^{12}=i\varphi_{M}$,  we get
\begin{equation}
\begin{array}{rcl}& &
\varphi_{n}^{12}(X_{1}(X_{2}W))= \sum\limits_{i+j=n}
\delta_{i}^{11}(X_{1})\varphi_{j}^{12}(X_{2}W)
\\\\&=&\sum\limits_{i+e=n} \delta_{i}^{11}(X_{1})\delta_{e}^{11}(X_{2})W +
\sum\limits_{i+e+k=n,k>0}
\delta_{i}^{11}(X_{1})\delta_{e}^{11}(X_{2})\varphi_{k}^{12}(W)
\end{array}
\end{equation}for any $X_{1}, X_{2}\in
{\mathcal{A}}$ and $W\in {\mathcal{M}}$. On the other hand
\begin{equation}
\begin{array}{rcl}& &
\varphi_{n}^{12}((X_{1}X_{2})W)= \sum\limits_{i+j=n,j>0}
\delta_{i}^{11}(X_{1}X_{2})\varphi_{j}^{12}(W)+\delta_{n}^{11}(X_{1}X_{2})W
\\\\&=&\sum\limits_{e+k+j=n,j>0}
\delta_{e}^{11}(X_{1})\delta_{k}^{11}(X_{2})\varphi_{j}^{12}(W)+\delta_{n}^{11}(X_{1}X_{2})W
\end{array}
\end{equation}for any $X_{1}, X_{2}\in
{\mathcal{A}}$ and $W\in {\mathcal{M}}$. Combining Eq. (12) with Eq.
(13), we get $[\delta_{n}^{11}(X_{1}X_{2})-\sum\limits_{e+i=n}
\delta_{i}^{11}(X_{1})\delta_{e}^{11}(X_{2})]W=0$. Since $M$ is
faithful, we get $\delta_{n}^{11}(X_{1}X_{2})=\sum\limits_{i+j=n}
\delta_{i}^{11}(X_{1})\delta_{j}^{11}(X_{2})$, i.e.
$\delta=\{\delta^{11}_{n}: n\in N\}$ is a higher derivation.

 Letting $S=\left [\begin{array}{cc}
  0 & -X^{-1}WY \\
  0 & Y \\
\end{array}\right ]$ and $T=\left
[\begin{array}{cc}
  X & W \\
  0 & 0 \\
\end{array}\right ]$ for any $Y\in {\mathcal{B}}$, $W\in {\mathcal{M}}$, and
invertible $X\in {\mathcal{A}}$. Then $ST=TS=0$. So we get
$$\begin{array}{rcl}& &\left [\begin{array}{cc}
  0 & 0 \\
  0 & 0 \\
\end{array}\right ]
=D_{n}(ST+TS)=\sum\limits_{i+j=n}(D_{i}(S)D_{j}(T)+D_{i}(T)D_{j}(S))
\\\\&=&\sum\limits_{i+j=n}(\left [\begin{array}{cc}
0 & -\varphi_{i}^{12}(X^{-1}WY)
+\tau_{i}^{12}(Y) \\
  0 & \tau_{i}^{22}(Y) \\
\end{array}\right ]\left [\begin{array}{ccc}
  \delta_{j}^{11}(X) & \delta_{j}^{12}(X)+\varphi_{j}^{12}(W) \\
  0 & 0 \\
\end{array}\right ]\\\\&&+\left [\begin{array}{ccc}
  \delta_{i}^{11}(X) &  \delta_{i}^{12}(X)+\varphi_{i}^{12}(W) \\
  0 & 0 \\
\end{array}\right ]\left [\begin{array}{cc}
 0 & -\varphi_{j}^{12}(X^{-1}WY)
+\tau_{j}^{12}(Y) \\
  0 & \tau_{j}^{22}(Y) \\
\end{array}\right ]).
\end{array}$$
The above equation implies that
\begin{equation}\
0=\sum\limits_{i+j=n}[\delta_{i}^{11}(X)(-\varphi_{j}^{12}(X^{-1}WY)+\tau_{j}^{12}(Y))+(\delta_{i}^{12}(X)+\varphi_{i}^{12}(W))\tau_{j}^{22}(Y)].
\end{equation}
By replacing $W$ by $\lambda W$ in the above equation, dividing the
equation by
 $\lambda$ and letting $\lambda\rightarrow +\infty$, we obtain that
\begin{equation}\
0=\sum\limits_{i+j=n}[-\delta_{i}^{11}(X)\varphi_{j}^{12}(X^{-1}WY)+\varphi_{i}^{12}(W)\tau_{j}^{22}(Y)].
\end{equation}
So we can get
\begin{equation}\
0=\sum\limits_{i+j=n}[-\delta_{i}^{11}(I_{1})\varphi_{j}^{12}(WY)+\varphi_{i}^{12}(W)\tau_{j}^{22}(Y)]
\end{equation}
by setting $X=I_{1}$ in the above equation. Since
$\delta=\{\delta^{11}_{n}: n\in N\}$ is a higher derivation,
$\delta_{n}^{11}(I_{1})=0$ when $n\geq 1$. It follows from Eq. (16)
that
\begin{equation}\
\varphi_{n}^{12}(WY)=\sum\limits_{i+j=n}\varphi_{i}^{12}(W)\tau_{j}^{22}(Y).
\end{equation}

We claim that $\tau=\{\tau^{22}_{n}: n\in N\}$ is a higher
derivation on $\mathcal{B}$. In fact, by the proof of [6, Theorem
2.1] we know that $\tau_{1}$ is a higher derivation. This implies
that
$\tau_{1}^{22}(Y_{1}Y_{2})=\tau_{1}^{22}(Y_{1})Y_{2}+Y_{1}\tau_{1}^{22}(Y_{2})$
for any $Y_{1}, Y_{2}\in {\mathcal{B}}$.  We now assume that
$\tau_{m}^{22}(Y_{1}Y_{2})= \sum\limits_{i+j=m}
\tau_{i}^{22}(Y_{1})\tau_{j}^{22}(Y_{2})$ for all $1\leq m<n$ with
$m\in N$. It follows from Eq. (17) that
\begin{equation}
\begin{array}{rcl}& &
\varphi_{n}^{12}(WY_{1}Y_{2})=
\varphi_{n}^{12}(W(Y_{1}Y_{2}))\\&=&W\tau_{n}^{22}(Y_{1}Y_{2})+\sum\limits_{i+j=n,j<n}
\varphi_{i}^{12}(W)\tau_{j}^{22}(Y_{1}Y_{2})
\\\\&=&W\tau_{n}^{22}(Y_{1}Y_{2})+\sum\limits_{i+e+k=n,i>0}
\varphi_{i}^{12}(W)\tau_{e}^{22}(Y_{1})\tau_{k}^{22}(Y_{2})
\end{array}
\end{equation} for any $Y_{1}, Y_{2}\in
{\mathcal{B}}$ and $W\in {\mathcal{M}}$. On the other hand by Eq.
(17) and the fact that $\mathcal{M}$ is a ($\mathcal{A}$,
$\mathcal{B}$)-bimodule, we have
\begin{equation}
\begin{array}{rcl}& &
\varphi_{n}^{12}(WY_{1}Y_{2})=
\varphi_{n}^{12}((WY_{1})Y_{2})\\\\&=&\sum\limits_{i+j=n}
\varphi_{i}^{12}(WY_{1})\tau_{j}^{22}(Y_{2}) =\sum\limits_{e+k+j=n}
\varphi_{e}^{12}(W)\tau_{k}^{22}(Y_{1})\tau_{j}^{22}(Y_{2})
\\\\&=&W\sum\limits_{k+j=n}
\tau_{e}^{22}(Y_{1})\tau_{j}^{22}(Y_{2})+\sum\limits_{e+k+j=n,e>0}
\varphi_{e}^{12}(W)\tau_{k}^{22}(Y_{1})\tau_{j}^{22}(Y_{2}).
\end{array}
\end{equation}
Combining Eq. (18) with Eq. (19), we get
$W[\tau_{n}^{22}(Y_{1}Y_{2})-\sum\limits_{k+j=n}
\tau_{e}^{22}(Y_{1})\tau_{j}^{22}(Y_{2})]W=0$. Since $M$ is
faithful, we get $\tau_{n}^{22}(Y_{1}Y_{2})=\sum\limits_{i+j=n}
\tau_{i}^{22}(Y_{1})\tau_{j}^{22}(Y_{2})$.

Now we prove that $(D_{n})_{n\in N}$ is a higher derivation. For any
$S=\left [\begin{array}{cc}
  X_{1} & W_{1}\\
  0 &  Y_{1} \\
\end{array}\right ], T=\left
[\begin{array}{cc}
  X_{2} & W_{2}\\
  0 &  Y_{2} \\
\end{array}\right ]\in \mathcal{T}$, where $X_{1},X_{2}\in
\mathcal{A}$, $W_{1},W_{2}\in \mathcal{M}$ and $Y_{1},Y_{2}\in
\mathcal{B}$. Summing up the above results and using the definition
of $D_{n}$, we obtain that
$$\begin{array}{rcl}D_{n}(ST)&=&D_{n}(\left [\begin{array}{cc}
  X_{1}X_{2} & X_{1}W_{2}+W_{1}Y_{2}\\
  0 & Y_{1}Y_{2} \\
\end{array}\right ])\\&=&\left [\begin{array}{cc}
  \delta_{n}^{11}(X_{1}X_{2}) & \delta_{n}^{12}(X_{1}X_{2})+\varphi_{n}^{12}(X_{1}W_{2}+W_{1}Y_{2})+\tau_{n}^{12}( Y_{1}Y_{2})  \\
  0 & \tau_{n}^{22}( Y_{1}Y_{2})  \\
\end{array}\right ],\end{array}$$
and$$\begin{array}{rcl}& &
\sum\limits_{i+j=n}D_{i}(S)D_{j}(T)=\sum\limits_{i+j=n}(\left
[\begin{array}{cc}
  \delta_{i}^{11}(X_{1}) & \delta_{i}^{12}(X_{1})+\varphi_{i}^{12}(W_{1})+\tau_{i}^{12}( Y_{1})  \\
  0 & \tau_{i}^{22}( Y_{1})  \\
\end{array}\right ]\\\\&&
\left [\begin{array}{cc}
 \delta_{j}^{11}(X_{2}) & \delta_{j}^{12}(X_{2})+\varphi_{j}^{12}(W_{2})+\tau_{j}^{12}( Y_{2})  \\
  0 & \tau_{j}^{22}( Y_{2})  \\
\end{array}\right ])
 \\\\&=&\left [\begin{array}{cc}
  \delta_{n}^{11}(X_{1}X_{2}) &
  \sum\limits_{i+j=n}(\delta_{i}^{11}(X_{1})\delta_{j}^{12}(X_{2})+
  \delta_{i}^{11}(X_{1})\tau_{j}^{12}( Y_{2})+\delta_{i}^{12}(X_{1})\tau_{j}^{22}(Y_{2})
  \\~ &+\tau_{i}^{12}(Y_{1})\tau_{j}^{22}( Y_{2}))+\varphi_{n}^{12}(X_{1}W_{2}+W_{1}Y_{2}) \\
   &\\
  0 & \tau_{n}^{22}( Y_{1}Y_{2})  \\
\end{array}\right ]
\end{array}$$
by Eq. (17) and the fact that both $\delta$ and $\tau$ are higher
derivations. So $D$ is a higher derivations if and only if the
equation
$$\begin{array}{rcl}& &
\delta_{n}^{12}(X_{1}X_{2})+\varphi_{n}^{12}(X_{1}W_{2}+W_{1}Y_{2})+\tau_{n}^{12}(
Y_{1}Y_{2})\\\\
&=&\sum\limits_{i+j=n}(\delta_{i}^{11}(X_{1})\delta_{j}^{12}(X_{2})+
  \delta_{i}^{11}(X_{1})\tau_{j}^{12}( Y_{2})\\&&+\delta_{i}^{12}(X_{1})\tau_{j}^{22}(Y_{2})
  +\tau_{i}^{12}(Y_{1})\tau_{j}^{22}( Y_{2}))+\varphi_{n}^{12}(X_{1}W_{2}+W_{1}Y_{2}) \\
\end{array}$$ holds.

We get that $\tau_{n}^{22}(I_{2})=0(n\geq1)$ from [4, lemma 2.2]. So
we can write
$$\delta_{n}^{12}(X)=-\sum\limits_{i+j=n}\delta_{i}^{11}(X)\tau_{j}^{12}(I_{2})$$
by setting $Y=I_{2}$ in Eq. (10). Letting $X=I_{1}$ in the above
equation, one gets $\delta_{n}^{12}(I_{1})=-\tau_{n}^{12}(I_{2})$.
So
\begin{equation}
\begin{array}{lll}
\delta_{n}^{12}(X)=
\sum\limits_{i+j=n}\delta_{i}^{11}(X)\delta_{j}^{12}(I_{1}).
\end{array}
\end{equation} Similarly by taking
$X=I_{1}$ in Eq. (10) and noting the fact $\delta_{n}^{11}(I_{1})=0
(n\geq 1)$, we have
\begin{equation}
\begin{array}{lll}
\tau_{n}^{12}(Y)=-\sum\limits_{i+j=n}\delta_{i}^{12}(I_{1})\tau_{j}^{22}(Y).
\end{array}
\end{equation}
Thus it follows from Eq. (20) and Eq. (21) that
\begin{equation}
\begin{array}{rcl}& &
\delta_{n}^{12}(X_{1}X_{2})+\tau_{n}^{12}(Y_{1}Y_{2})=
\sum\limits_{i+j=n}\delta_{i}^{11}(X_{1}X_{2})\delta_{j}^{12}(I_{1})-
\sum\limits_{i+j=n}\delta_{i}^{12}(I_{1})\tau_{j}^{22}(Y_{1}Y_{2})
\\\\&=&\sum\limits_{k+l+j=n}\delta_{k}^{11}(X_{1})\delta_{l}^{11}(X_{2})\delta_{j}^{12}(I_{1})-
\sum\limits_{i+k+l=n}\delta_{i}^{12}(I_{1})\tau_{k}^{22}(Y_{1})\tau_{l}^{22}(Y_{2}).
\end{array}
\end{equation}
On the other hand
\begin{equation}
\begin{array}{rcl}& &
\sum\limits_{i+j=n}(\delta_{i}^{11}(X_{1})\delta_{j}^{12}(X_{2})
+\delta_{i}^{11}(X_{1})\tau_{j}^{12}(Y_{2})
+\delta_{i}^{12}(X_{1})\tau_{j}^{22}(Y_{2})
+\tau_{i}^{12}(Y_{1})\tau_{j}^{22}(Y_{2}))\\\\
&=&\sum\limits_{i+j=n}\sum\limits_{k+l=j}\delta_{i}^{11}(X_{1})\delta_{k}^{11}(X_{2})\delta_{l}^{12}(I_{1})-
\sum\limits_{i+j=n}\sum\limits_{k+l=j}\delta_{i}^{11}(X_{1})\delta_{k}^{12}(I_{1})\tau_{l}^{22}(Y_{2})
\\\\&&+\sum\limits_{i+j=n}\sum\limits_{k+l=i}\delta_{k}^{11}(X_{1})\delta_{l}^{12}(I_{1})\tau_{j}^{22}(Y_{2})-
\sum\limits_{i+j=n}\sum\limits_{k+l=i}\delta_{k}^{11}(I_{1})\tau_{l}^{22}(Y_{1})\tau_{j}^{22}(Y_{2})
\\\\&=&\sum\limits_{i+k+l=n}\delta_{i}^{11}(X_{k})\delta_{k}^{11}(X_{2})\delta_{l}^{12}(I_{1})-
\sum\limits_{j+k+l=n}\delta_{k}^{12}(I_{1})\tau_{l}^{22}(Y_{1})\tau_{j}^{22}(Y_{2}).
\end{array}
\end{equation}
Thus combining Eq. (22) with Eq. (23), we arrive at
$$\begin{array}{rcl}& &
\delta_{n}^{12}(X_{1}X_{2})+\varphi_{n}^{12}(X_{1}W_{2}+W_{1}Y_{2})+\tau_{n}^{12}(
Y_{1}Y_{2})\\\\
&=&\sum\limits_{i+j=n}(\delta_{i}^{11}(X_{1})\delta_{j}^{12}(X_{2})+
  \delta_{i}^{11}(X_{1})\tau_{j}^{12}( Y_{2})+\delta_{i}^{12}(X_{1})\tau_{j}^{22}(Y_{2})
  \\&&+\tau_{i}^{12}(Y_{1})\tau_{j}^{22}( Y_{2}))+\varphi_{n}^{12}(X_{1}W_{2}+W_{1}Y_{2}). \\
\end{array}$$ Finally we obtain the desired result.\\~\\
{\bf{Theorem 2.2}} {{\it Let $D=\{D_{n}\}$ be a family of additive
mappings on $\mathcal{T}$ that $D_{0}=iD_{\mathcal{T}}$. If $D$ is
Jordan higher derivable at $G=\left [\begin{array}{cc}
  I_{1} & X_{0} \\
  0 & I_{2} \\
\end{array}\right ]$, then $D$ is a higher derivation.}\\
{\bf Proof.} We set $S=\left [\begin{array}{cc}
  X & 0 \\
  0 & Y \\
\end{array}\right ]$ and $T=\left
[\begin{array}{cc}
  X^{-1} & X^{-1}X_{0} \\
  0 & Y^{-1}\\
\end{array}\right ]$ for every invertible element $X\in {\mathcal{A}}$ and
$Y\in {\mathcal{B}}$. Then $ST=G$ and $TS=\left [\begin{array}{cc}
  I_{1} & X^{-1}X_{0}Y \\
  0 & I_{2}\\
\end{array}\right ]$, so we obtain
$$\begin{array}{rcl}& &\left
[\begin{array}{cc}
  2\delta_{n}^{11}(I_{1})+2\tau_{n}^{11}(I_{2}) & 2\delta_{n}^{12}(I_{1})+2\tau_{n}^{12}(I_{2})+
  \\
+\varphi_{n}^{11}(X_{0}+X^{-1}X_{0}Y)&+\varphi_{n}^{12}(X_{0}+X^{-1}X_{0}Y)
\\
~&~\\
  0 & 2\delta_{n}^{22}(I_{1})+\varphi_{n}^{22}(X_{0}+X^{-1}X_{0}Y)+2\tau_{n}^{22}(I_{2})\\
\end{array}\right ]\\\\&=&D_{n}(ST+TS)=\sum\limits_{i+j=n}(D_{i}(S)D_{j}(T)+D_{i}(T)D_{j}(S))\\
\\&=&\sum\limits_{i+j=n}(\left [\begin{array}{cc}
  \delta_{i}^{11}(X)+\tau_{i}^{11}(Y) & \delta_{i}^{12}(X)+\tau_{i}^{12}(Y) \\
  0 & \delta_{i}^{22}(X)+\tau_{i}^{22}(Y) \\
\end{array}\right ]
\\\\&&\left [\begin{array}{cc}
  \delta_{j}^{11}(X^{-1})+\varphi_{j}^{11}(X^{-1}X_{0}) & \delta_{j}^{12}(X^{-1})+\varphi_{j}^{12}(X^{-1}X_{0}) \\
  +\tau_{j}^{11}(Y^{-1})&+\tau_{j}^{12}(Y^{-1})\\\\
  0 & \delta_{j}^{22}(X^{-1})+\varphi_{j}^{22}(X^{-1}X_{0})+\tau_{j}^{22}(Y^{-1}) \\
\end{array}\right ]\\\\&&+\left [\begin{array}{cc}
 \delta_{i}^{11}(X^{-1})+\varphi_{i}^{11}(X^{-1}X_{0}) & \delta_{i}^{12}(X^{-1})+\varphi_{i}^{12}(X^{-1}X_{0}) \\
+\tau_{i}^{11}(Y^{-1})&+\tau_{i}^{12}(Y^{-1})\\\\
  0 & \delta_{i}^{22}(X^{-1})+\varphi_{i}^{22}(X^{-1}X_{0})+\tau_{i}^{22}(Y^{-1}) \\
\end{array}\right ]
\\\\&&\left [\begin{array}{cc}
  \delta_{j}^{11}(X)+\tau_{j}^{11}(Y) & \delta_{j}^{12}(X)+\tau_{j}^{12}(Y) \\
  0 & \delta_{j}^{22}(X)+\tau_{j}^{22}(Y) \\
\end{array}\right ]).\\
\end{array}$$
So according to the above matrix equation, we get
\begin{equation}\label{10}\begin{array}{rcl}& &
2\delta_{n}^{11}(I_{1})+2\tau_{n}^{11}(I_{2})+\varphi_{n}^{11}(X_{0}+X^{-1}X_{0}Y)\\\\
&=&\sum\limits_{i+j=n}[(\delta_{i}^{11}(X)+\tau_{i}^{11}(Y))
(\delta_{j}^{11}(X^{-1})+\varphi_{j}^{11}(X^{-1}X_{0})+\tau_{j}^{11}(Y^{-1})
 )\\\\&&
+(\delta_{i}^{11}(X^{-1})+\varphi_{i}^{11}(X^{-1}X_{0})+\tau_{i}^{11}(Y^{-1}))
(\delta_{j}^{11}(X)+\tau_{j}^{11}(Y))],
\end{array}
\end{equation}

\begin{equation}\label{10}\begin{array}{rcl}& &
2\delta_{n}^{12}(I_{1})+2\tau_{n}^{12}(I_{2})+\varphi_{n}^{12}(X_{0}+X^{-1}X_{0}Y)\\\\
&=&\sum\limits_{i+j=n}[(\delta_{i}^{11}(X)+\tau_{i}^{11}(Y))
(\delta_{j}^{12}(X^{-1})+\varphi_{j}^{12}(X^{-1}X_{0})+\tau_{j}^{12}(Y^{-1})
 )\\\\&&
+(\delta_{i}^{12}(X)+\tau_{i}^{12}(Y))
(\delta_{j}^{22}(X^{-1})+\varphi_{j}^{22}(X^{-1}X_{0})+\tau_{j}^{22}(Y^{-1})
 )\\\\&&
+(\delta_{i}^{11}(X^{-1})+\varphi_{i}^{11}(X^{-1}X_{0})+\tau_{i}^{11}(Y^{-1}))
(\delta_{j}^{12}(X)+\tau_{j}^{12}(Y))\\\\&&
+(\delta_{i}^{12}(X^{-1})+\varphi_{i}^{12}(X^{-1}X_{0})+\tau_{i}^{12}(Y^{-1}))
(\delta_{j}^{22}(X)+\tau_{j}^{22}(Y))],
\end{array}
\end{equation}\\

\begin{equation}\label{10}\begin{array}{rcl}& &
2\delta_{n}^{22}(I_{1})+2\tau_{n}^{22}(I_{2})+\varphi_{n}^{22}(X_{0}+X^{-1}X_{0}Y)\\\\
&=&\sum\limits_{i+j=n}[(\delta_{i}^{22}(X)+\tau_{i}^{22}(Y))
(\delta_{j}^{22}(X^{-1})+\varphi_{j}^{22}(X^{-1}X_{0})+\tau_{j}^{22}(Y^{-1})
 )\\\\&&
+(\delta_{i}^{22}(X^{-1})+\varphi_{i}^{22}(X^{-1}X_{0})+\tau_{i}^{22}(Y^{-1}))
(\delta_{j}^{22}(X)+\tau_{j}^{22}(Y))].
\end{array}
\end{equation}

We claim that
$\delta_{n}^{11}(I_{1})=\tau_{n}^{11}(I_{2})=\varphi_{n}^{11}(X_{0})=0$
when $n\geq1$ . In fact, we could obtain
\begin{equation}\label{10}\begin{array}{rcl}& &
2\delta_{n}^{11}(I_{1})+2\tau_{n}^{11}(I_{2})+\varphi_{n}^{11}(X_{0}+X_{0})\\\\
&=&\sum\limits_{i+j=n}[(\delta_{i}^{11}(I_{1})+\tau_{i}^{11}(I_{2}))
(\delta_{j}^{11}(I_{1})+\varphi_{j}^{11}(X_{0})+\tau_{j}^{11}(I_{2})
 )\\\\&&
+(\delta_{i}^{11}(I_{1})+\varphi_{i}^{11}(X_{0})+\tau_{i}^{11}(I_{2}))
(\delta_{j}^{11}(I_{1})+\tau_{j}^{11}(I_{2}))]
\end{array}
\end{equation}
by setting $X=I_{1}$ and $Y=I_{2}$ in Eq. (24). When $n=1$, the
result that
$\delta_{1}^{11}(I_{1})=\tau_{1}^{11}(I_{2})=\varphi_{1}^{11}(X_{0})=0$
holds according to the [6, Theorem 2.2]. So we assume that
$\delta_{m}^{11}(I_{1})=\tau_{m}^{11}(I_{2})=\varphi_{m}^{11}(X_{0})=0$
for all $1\leq m <n, m\in N$. Combining Eq. (27)  with the fact
$\delta_{0}^{11}(I_{1})=I_{1},\tau_{0}^{11}(I_{2})=0$ and using the
induction hypothesis, we have
$$\begin{array}{rcl}& &
2\delta_{n}^{11}(I_{1})+2\tau_{n}^{11}(I_{2})+2\varphi_{n}^{11}(X_{0})
=\delta_{n}^{11}(I_{1})+\tau_{n}^{11}(I_{2})+\delta_{n}^{11}(I_{1})+\tau_{n}^{11}(I_{2})\\\\&&
+2\delta_{n}^{11}(I_{1})+2\tau_{n}^{11}(I_{2})+2\varphi_{n}^{11}(X_{0}).
\end{array}
$$
Hence $\delta_{n}^{11}(I_{1})+\tau_{n}^{11}(I_{2})=0(n\geq 1)$.
Similarly we also can set that $X=I_{1}$ and $Y=-I_{2}$ in Eq. (24).
Using the induction hypothesis, we get
$\delta_{n}^{11}(I_{1})-\tau_{n}^{11}(I_{2})=-\varphi_{n}^{11}(X_{0})$.
Summing up the above equations we get
$2\delta_{n}^{11}(I_{1})=-2\tau_{n}^{11}(I_{2})=\varphi_{n}^{11}(X_{0})$.

Setting $X=\frac{1}{2}I_{1}$ and $Y=I_{2}$ in Eq. (24) and using
$\delta_{n}^{11}(I_{1})+\tau_{n}^{11}(I_{2})=0$, we have
$$\begin{array}{rcl}& &
3\varphi_{n}^{11}(X_{0})=\sum\limits_{i+j=n}[(\frac{1}{2}\delta_{i}^{11}(I_{1})+\tau_{i}^{11}(I_{2}))
(2\delta_{j}^{11}(I_{1})+\tau_{j}^{11}(I_{2})+2\varphi_{j}^{11}(X_{0})
)\\\\
&&+(2\delta_{i}^{11}(I_{1})+\tau_{i}^{11}(I_{2})+2\varphi_{i}^{11}(X_{0})
)(\frac{1}{2}\delta_{j}^{11}(I_{1})+\tau_{j}^{11}(I_{2}))].
\end{array}$$

Thus combining
$2\delta_{n}^{11}(I_{1})=-2\tau_{n}^{11}(I_{2})=\varphi_{n}^{11}(X_{0})$
with the assumption and using $\delta_{0}^{11}(I_{1})=I_{1}$, one
obtains
$$\begin{array}{rcl}& &
3\varphi_{n}^{11}(X_{0})=\frac{1}{2}(2\delta_{n}^{11}(I_{1})+\tau_{n}^{11}(I_{2})+2\varphi_{n}^{11}(X_{0}))\\\\
&&+2(\delta_{n}^{11}(I_{1})+\tau_{n}^{11}(I_{2}))+2(\delta_{n}^{11}(I_{1})+\tau_{n}^{11}(I_{2}))\\\\
&&+\frac{1}{2}(2\delta_{n}^{11}(I_{1})+\tau_{n}^{11}(I_{2})+2\varphi_{n}^{11}(X_{0}))
\end{array}.$$
So
$\varphi_{n}^{11}(X_{0})=4\delta_{n}^{11}(I_{1})+5\tau_{n}^{11}(I_{2})$.
We can claim that
$\delta_{n}^{11}(I_{1})=\tau_{n}^{11}(I_{2})=\varphi_{n}^{11}(X_{0})=0$.
Hence the Eq. (24) can be rewritten into
\begin{equation}\label{10}\begin{array}{rcl}& &
\varphi_{n}^{11}(X^{-1}X_{0}Y)
=\sum\limits_{i+j=n}[(\delta_{i}^{11}(X)+\tau_{i}^{11}(Y))
(\varphi_{j}^{11}(X^{-1}X_{0})+\tau_{j}^{11}(Y^{-1})+
\delta_{j}^{11}(X^{-1}) )\\\\&&
+(\delta_{i}^{11}(X^{-1})+\varphi_{i}^{11}(X^{-1}X_{0})+\tau_{i}^{11}(Y^{-1}))
(\delta_{j}^{11}(X)+\tau_{j}^{11}(Y))].
\end{array}
\end{equation}

Similarly by setting $X=I_{1}$ and $Y=I_{2}$ in Eq. (26) and using
the induction, we can get
$\delta_{n}^{22}(I_{1})+\tau_{n}^{22}(I_{2})=0$. We also can obtain
$\delta_{n}^{22}(I_{1})=\tau_{n}^{22}(I_{2})=\varphi_{n}^{22}(X_{0})=0$
if we take $X=I_{1}$ and $Y=\frac{1}{2}I_{2}$ in Eq. (27). Thus
\begin{equation}\label{10}\begin{array}{rcl}& &
\varphi_{n}^{22}(X^{-1}X_{0}Y)
=\sum\limits_{i+j=n}[(\delta_{i}^{22}(X)+\tau_{i}^{22}(Y))
(\varphi_{j}^{22}(X^{-1}X_{0})+\tau_{j}^{22}(Y^{-1})+
\delta_{j}^{22}(X^{-1}) )\\\\&&
+(\delta_{i}^{22}(X^{-1})+\varphi_{i}^{22}(X^{-1}X_{0})+\tau_{i}^{22}(Y^{-1}))
(\delta_{j}^{22}(X)+\tau_{j}^{22}(Y))].
\end{array}
\end{equation}

We take $X=I_{1}$ and $Y=I_{2}$ in Eq. (25), then we can get
$\delta_{n}^{12}(I_{1})+\tau_{n}^{12}(I_{2})=0$. Letting
respectively $Y=I_{2}$ and $Y=\frac{1}{2}I_{2}$ in Eq. (25) and
using the above equation we have

\begin{equation}\label{10}\begin{array}{rcl}& &
\varphi_{n}^{12}(X_{0}+X^{-1}X_{0})
=\sum\limits_{i+j=n}[\delta_{i}^{11}(X)
(\delta_{j}^{12}(X^{-1})+\varphi_{j}^{12}(X^{-1}X_{0})+
\tau_{j}^{12}(I_{2}) )\\\\&&
+(\delta_{i}^{12}(X)+\tau_{i}^{12}(I_{2})) (\delta_{j}^{22}(X^{-1})+
\varphi_{j}^{22}(X^{-1}X_{0}) )\\\\&&
+(\delta_{i}^{11}(X^{-1})+\varphi_{i}^{11}(X^{-1}X_{0}))
(\delta_{j}^{12}(X)+\tau_{j}^{12}(I_{2}))\\\\&&
+(\delta_{i}^{12}(X^{-1})+\varphi_{i}^{12}(X^{-1}X_{0})+\tau_{i}^{12}(I_{2}))
\delta_{j}^{22}(X)]\\\\&& +\delta_{n}^{12}(X)+\tau_{n}^{12}(I_{2})+
\delta_{n}^{12}(X^{-1})+\varphi_{n}^{12}(X^{-1}X_{0})+\tau_{n}^{12}(I_{2}),
\end{array}
\end{equation}

\begin{equation}\label{10}\begin{array}{rcl}& &
\varphi_{n}^{12}(X_{0}+\frac{1}{2}X^{-1}X_{0})
=\sum\limits_{i+j=n}[\delta_{i}^{11}(X) (\delta_{j}^{12}(X^{-1})
+\varphi_{j}^{12}(X^{-1}X_{0})+2\tau_{j}^{12}(I_{2}))\\\\&&
+(\delta_{i}^{12}(X)+\frac{1}{2}\tau_{i}^{12}(I_{2}))
(\delta_{j}^{22}(X^{-1})+ \varphi_{j}^{22}(X^{-1}X_{0}))\\\\&&
+(\delta_{i}^{11}(X^{-1})+\varphi_{i}^{11}(X^{-1}X_{0}))
(\delta_{j}^{12}(X)+\frac{1}{2}\tau_{j}^{12}(I_{2}))\\\\&&
+(\delta_{i}^{12}(X^{-1})+\varphi_{i}^{12}(X^{-1}X_{0})+2\tau_{i}^{12}(I_{2}))
\delta_{j}^{22}(X)]\\\\&& +2\delta_{n}^{12}(X)+\tau_{n}^{12}(I_{2})+
\frac{1}{2}\delta_{n}^{12}(X^{-1})+\frac{1}{2}\varphi_{n}^{12}(X^{-1}X_{0})+\tau_{n}^{12}(I_{2}),
\end{array}
\end{equation}
which implies that
$$\begin{array}{rcl}& &
\frac{1}{2}\varphi_{n}^{12}(X^{-1}X_{0})
=\sum\limits_{i+j=n}[-\delta_{i}^{11}(X)\tau_{j}^{12}(I_{2})\\\\
&&+\frac{1}{2}\tau_{i}^{12}(I_{2})(\delta_{j}^{22}(X^{-1})+\varphi_{j}^{22}(X^{-1}X_{0}))+\frac{1}{2}(\delta_{i}^{11}(X^{-1})+\varphi_{i}^{11}(X^{-1}X_{0}))\tau_{j}^{12}(I_{2})\\\\
&&-\tau_{i}^{12}(I_{2})\delta_{j}^{22}(X)]-\delta_{n}^{12}(X)+\frac{1}{2}\delta_{n}^{12}(X^{-1})+\frac{1}{2}\varphi_{n}^{12}(X^{-1}X_{0}).
\end{array}$$
So
\begin{equation}\label{10}\begin{array}{rcl}& &
\frac{1}{2}\sum\limits_{i+j=n}[\tau_{i}^{12}(I_{2})\delta_{j}^{22}(X^{-1})
+\delta_{i}^{11}(X^{-1})\tau_{j}^{12}(I_{2})\\\\&&+\tau_{i}^{12}(I_{2})\varphi_{j}^{22}(X^{-1}X_{0})
+\varphi_{i}^{11}(X^{-1}X_{0})\tau_{j}^{12}(I_{2})]+\frac{1}{2}\delta_{n}^{12}(X^{-1})\\\\
&=&\sum\limits_{i+j=n}[\delta_{i}^{11}(X)\tau_{j}^{12}(I_{2})+\tau_{i}^{12}(I_{2})\delta_{j}^{22}(X)]
+\delta_{n}^{12}(X).
\end{array}
\end{equation}
Thus we get
\begin{equation}\label{10}\begin{array}{rcl}& &
\frac{1}{2}\sum\limits_{i+j=n}[\tau_{i}^{12}(I_{2})\delta_{j}^{22}(X)
+\delta_{i}^{11}(X)\tau_{j}^{12}(I_{2})\\\\&&+\tau_{i}^{12}(I_{2})\varphi_{j}^{22}(XX_{0})
+\varphi_{i}^{11}(XX_{0})\tau_{j}^{12}(I_{2})]+\frac{1}{2}\delta_{n}^{12}(X)\\\\
&=&\sum\limits_{i+j=n}[\delta_{i}^{11}(X^{-1})\tau_{j}^{12}(I_{2})+\tau_{i}^{12}(I_{2})\delta_{j}^{22}(X^{-1})]
+\delta_{n}^{12}(X^{-1})
\end{array}
\end{equation}
for any invertible $X\in \mathcal{A}$ by replacing $X^{-1}$ by $X$
in Eq.(32). It follows that
$$\begin{array}{rcl}& &
\frac{1}{2}[\frac{1}{2}\sum\limits_{i+j=n}[\tau_{i}^{12}(I_{2})\delta_{j}^{22}(X)
+\delta_{i}^{11}(X)\tau_{j}^{12}(I_{2})\\\\&&+\tau_{i}^{12}(I_{2})\varphi_{j}^{22}(XX_{0})
+\varphi_{i}^{11}(XX_{0})\tau_{j}^{12}(I_{2})]+\frac{1}{2}\delta_{n}^{12}(X)]\\\\
&&+\frac{1}{2}\sum\limits_{i+j=n}[\tau_{i}^{12}(I_{2})\varphi_{j}^{22}(X^{-1}X_{0})
+\varphi_{i}^{11}(X^{-1}X_{0})\tau_{j}^{12}(I_{2})]\\\\
&&=\sum\limits_{i+j=n}[\delta_{i}^{11}(X)\tau_{j}^{12}(I_{2})+\tau_{i}^{12}(I_{2})\delta_{j}^{22}(X)]
+\delta_{n}^{12}(X).
\end{array}$$
So
\begin{equation}\label{10}\begin{array}{rcl}& &
\frac{1}{4}[\sum\limits_{i+j=n}[\tau_{i}^{12}(I_{2})\delta_{j}^{22}(X)
+\delta_{i}^{11}(X)\tau_{j}^{12}(I_{2})]+\delta_{n}^{12}(X)]\\\\&&+
\frac{1}{4}\sum\limits_{i+j=n}[\tau_{i}^{12}(I_{2})\varphi_{j}^{22}(XX_{0})
+\varphi_{i}^{11}(XX_{0})\tau_{j}^{12}(I_{2})]
\\\\&&+\frac{1}{2}\sum\limits_{i+j=n}[\tau_{i}^{12}(I_{2})\varphi_{j}^{22}(X^{-1}X_{0})
+\varphi_{i}^{11}(X^{-1}X_{0})\tau_{j}^{12}(I_{2})]
\\\\&=&\sum\limits_{i+j=n}[\tau_{i}^{12}(I_{2})\delta_{j}^{22}(X)
+\delta_{i}^{11}(X)\tau_{j}^{12}(I_{2})]+\delta_{n}^{12}(X)
\end{array}
\end{equation}
for any invertible $X\in \mathcal{A}$.

Similarly by letting  $X=I_{1}$ and $X=2I_{1}$ in Eq. (25), it is
easily checked that
\begin{equation}\label{10}\begin{array}{rcl}& &
\varphi_{n}^{12}(X_{0}+X_{0}Y) =\sum\limits_{i+j=n}[\tau_{i}^{11}(Y)
(\varphi_{j}^{12}(X_{0})+\tau_{j}^{12}(Y^{-1})+
\delta_{j}^{12}(I_{1}))
\\\\&&+(\delta_{i}^{12}(I_{1})+\tau_{i}^{12}(Y))
\tau_{j}^{22}(Y^{-1})+\tau_{i}^{11}(Y^{-1})(\delta_{j}^{12}(I_{1})+\tau_{j}^{12}(Y))
\\\\&&+(\varphi_{i}^{12}(X_{0})+\tau_{i}^{12}(Y^{-1})+
\delta_{i}^{12}(I_{1}))\tau_{j}^{22}(Y)]
\\\\&&+\varphi_{n}^{12}(X_{0})+\tau_{n}^{12}(Y^{-1})+
2\delta_{n}^{12}(I_{1})+\tau_{n}^{12}(Y),
\end{array}
\end{equation}\\

\begin{equation}\label{10}\begin{array}{rcl}& &
\varphi_{n}^{12}(X_{0}+\frac{1}{2}X_{0}Y)
=\sum\limits_{i+j=n}[\tau_{i}^{11}(Y)
(\frac{1}{2}\varphi_{j}^{12}(X_{0})+\tau_{j}^{12}(Y^{-1})+
\frac{1}{2}\delta_{j}^{12}(I_{1}))
\\\\&&+(2\delta_{i}^{12}(I_{1})+\tau_{i}^{12}(Y))
\tau_{j}^{22}(Y^{-1})+\tau_{i}^{11}(Y^{-1})(2\delta_{j}^{12}(I_{1})+\tau_{j}^{12}(Y))
\\\\&&+(\frac{1}{2}\varphi_{i}^{12}(X_{0})+\tau_{i}^{12}(Y^{-1})+
\frac{1}{2}\delta_{i}^{12}(I_{1}))\tau_{j}^{22}(Y)]
\\\\&&+\varphi_{n}^{12}(X_{0})+2\tau_{n}^{12}(Y^{-1})+
2\delta_{n}^{12}(I_{1})+\frac{1}{2}\tau_{n}^{12}(Y),
\end{array}
\end{equation}
which implies that
\begin{equation}\label{10}\begin{array}{rcl}& &
\frac{1}{2}\varphi_{n}^{12}(X_{0}Y)=
\sum\limits_{i+j=n}[\frac{1}{2}\tau_{i}^{11}(Y)
(\varphi_{j}^{12}(X_{0})+\delta_{j}^{12}(I_{1}))
-\delta_{i}^{12}(I_{1})\tau_{j}^{22}(Y^{-1})
\\\\&&-\tau_{i}^{11}(Y^{-1})\delta_{j}^{12}(I_{1})
+\frac{1}{2}(\varphi_{i}^{12}(X_{0})+\delta_{i}^{12}(I_{1}))\tau_{j}^{22}(Y)]
+\frac{1}{2}\tau_{n}^{12}(Y)-\tau_{n}^{12}(Y^{-1}).
\end{array}
\end{equation}

By considering Eq. (28) and $\varphi_{n}^{11}(X_{0})=0$ and letting
$X=I_{1}$ and $X=2I_{1}$ respectively, it is easily verified that
\begin{equation}\label{10}\begin{array}{rcl}& &
\varphi_{n}^{11}(X_{0}Y)=\sum\limits_{i+j=n}[\tau_{i}^{11}(Y)\tau_{j}^{11}(Y^{-1})
+\tau_{i}^{11}(Y^{-1})\tau_{j}^{11}(Y)]+2\tau_{n}^{11}(Y^{-1})+2\tau_{n}^{11}(Y),
\end{array}
\end{equation}

\begin{equation}\label{10}\begin{array}{rcl}& &
\frac{1}{2}\varphi_{n}^{11}(X_{0}Y)=\sum\limits_{i+j=n}[\tau_{i}^{11}(Y)\tau_{j}^{11}(Y^{-1})
+\tau_{i}^{11}(Y^{-1})\tau_{j}^{11}(Y)]+4\tau_{n}^{11}(Y^{-1})+\tau_{n}^{11}(Y).
\end{array}
\end{equation}
When $n=0$, $\tau_{0}^{11}(Y)=0$. When $n=1$, $\tau_{1}^{11}(Y)=0$
according to [6, Theorem 2.2]. We assume that $\tau_{m}^{11}(Y)=0$
for any $Y\in \mathcal{B}$ and $1\leq m<n$. So combining Eq. (38)
with Eq. (39) and using the induction hypothesis, we have
\begin{equation}
\varphi_{n}^{11}(X_{0}Y)=2\tau_{n}^{11}(Y^{-1})+2\tau_{n}^{11}(Y),
\end{equation}
\begin{equation}
\frac{1}{2}\varphi_{n}^{11}(X_{0}Y)=4\tau_{n}^{11}(Y^{-1})+\tau_{n}^{11}(Y).
\end{equation}
By direct computation, one can verify that
$\tau_{n}^{11}(Y^{-1})=0$. There exists $n\in N$ such that
$nI_{2}-Y$ is invertible for any $Y\in\mathcal{B}$ and
$\tau_{n}^{11}(I_{2})=0$, so $\tau_{n}^{11}(Y)=0$  for any
$Y\in\mathcal{B}$ .

When $n=0$, $\delta_{0}^{22}(X)=0$ for any $X\in \mathcal{A}$. By
[6, Theorem 2.2], we can claim that  When $n=1$,
$\delta_{1}^{22}(X)=0$. So now we assume that $\delta_{m}^{22}(X)=0$
for all $1\leq m<n$ and $X\in \mathcal{A}$. Taking respectively
$Y=I_{2}$ and $Y=2I_{2}$ in Eq. (29) and using
$\tau_{n}^{22}(I_{2})=0, n\geq 1$,
$\tau_{0}^{22}=i\tau_{\mathcal{B}}$ we have
\begin{equation}\label{10}\begin{array}{rcl}& &
\varphi_{n}^{22}(X^{-1}X_{0})=\sum\limits_{i+j=n}
[\delta_{i}^{22}(X)(\delta_{j}^{22}(X^{-1})+\varphi_{j}^{22}(X^{-1}X_{0}))
\\\\&&+(\varphi_{i}^{22}(X^{-1}X_{0})+\delta_{i}^{22}(X^{-1}))\delta_{j}^{22}(X)]
\\\\&&+2\delta_{n}^{22}(X)+2\varphi_{n}^{22}(X^{-1}X_{0})+2\delta_{n}^{22}(X^{-1}),
\end{array}
\end{equation}
and

\begin{equation}\label{10}\begin{array}{rcl}& &
2\varphi_{n}^{22}(X^{-1}X_{0})=\sum\limits_{i+j=n}
[\delta_{i}^{22}(X)(\delta_{j}^{22}(X^{-1})+\varphi_{j}^{22}(X^{-1}X_{0}))
\\\\&&+(\varphi_{i}^{22}(X^{-1}X_{0})+\delta_{i}^{22}(X^{-1}))\delta_{j}^{22}(X)]
\\\\&&+\delta_{n}^{22}(X)+4\varphi_{n}^{22}(X^{-1}X_{0})+4\delta_{n}^{22}(X^{-1}).
\end{array}
\end{equation}
Combining the assumption and the above equations, we have the
following equations:
$$-\varphi_{n}^{22}(X^{-1}X_{0})=2\delta_{n}^{22}(X)+2\delta_{n}^{22}(X^{-1}),$$
$$-2\varphi_{n}^{22}(X^{-1}X_{0})=\delta_{n}^{22}(X)+4\delta_{n}^{22}(X^{-1}).$$
By direct computation, one can verify that $\delta_{n}^{22}(X)=0$
for any invertible $X\in\mathcal{A}$ and $n\in N$. Because there is
some integer $n$ such that $nI_{1}-X$ is invertible for every
$X\in\mathcal{A}$, the conclusion of $\delta_{n}^{22}(X)=0$ holds
for every $X\in\mathcal{A}$.

We set $S=\left [\begin{array}{cc}
  X & XW \\
  0 & Y \\
\end{array}\right ]$ and
$T=\left[\begin{array}{cc}
  X^{-1} & X^{-1}X_{0}-WY^{-1} \\
  0 & Y^{-1}\\
\end{array}\right ]$ for any $Y\in {\mathcal{B}}$, $W\in {\mathcal{M}}$, and for
any invertible $X\in {\mathcal{A}}$, then $ST=G$ and
$TS=\left[\begin{array}{cc}
  I_{1} & X^{-1}X_{0}Y \\
  0 & I_{2}\\
\end{array}\right ]$. So combining
$\delta_{n}^{12}(I_{1})+\tau_{n}^{12}(I_{2})=0$ with the
characterization of $D$,  we obtain the following when $n\geq1$

$$\begin{array}{rcl}& &
\left[\begin{array}{cc}
 \varphi_{n}^{11}(X^{-1}X_{0}Y) & \varphi_{n}^{12}(X_{0}+X^{-1}X_{0}Y) \\
 0 & \varphi_{n}^{22}(X^{-1}X_{0}Y)\\
\end{array}\right ]\\\\&=&D_{n}(ST+TS)=\sum\limits_{i+j=n}(D_{i}(S)D_{j}(T)+D_{i}(T)D_{j}(S))\\
\\&=&\sum\limits_{i+j=n}(\left [\begin{array}{cc}
  \delta_{i}^{11}(X)+\varphi_{i}^{11}(XW) & \delta_{i}^{12}(X)+\varphi_{i}^{12}(XW)+\tau_{i}^{12}(Y) \\
  0 & \tau_{i}^{22}(Y)+\varphi_{i}^{22}(XW) \\
\end{array}\right ]
\\\\&&\left [\begin{array}{cc}
  \delta_{j}^{11}(X^{-1})+\varphi_{j}^{11}(X^{-1}X_{0}-WY^{-1}) & \delta_{j}^{12}(X^{-1})+\varphi_{j}^{12}(X^{-1}X_{0}-WY^{-1})+\tau_{j}^{12}(Y) \\
  0 & \tau_{j}^{22}(Y^{-1})+\varphi_{j}^{22}(X^{-1}X_{0}-WY^{-1}) \\
\end{array}\right ]\\\\&&+\left [\begin{array}{cc}
 \delta_{i}^{11}(X^{-1})+\varphi_{i}^{11}(X^{-1}X_{0}-WY^{-1}) & \delta_{i}^{12}(X^{-1})+\varphi_{i}^{12}(X^{-1}X_{0}-WY^{-1})+\tau_{i}^{12}(Y) \\
  0 & \tau_{i}^{22}(Y^{-1})+\varphi_{i}^{22}(X^{-1}X_{0}-WY^{-1}) \\
\end{array}\right ]
\\\\&&\left [\begin{array}{cc}
   \delta_{j}^{11}(X)+\varphi_{j}^{11}(XW) & \delta_{j}^{12}(X)+\varphi_{j}^{12}(XW)+\tau_{j}^{12}(Y) \\
  0 & \tau_{j}^{22}(Y)+\varphi_{j}^{22}(XW) \\
\end{array}\right ],\\
\end{array}$$
which implies the following three equations
\begin{equation}\label{10}\begin{array}{rcl}& &
\varphi_{n}^{11}(X^{-1}X_{0}Y)=\sum\limits_{i+j=n}[(\delta_{i}^{11}(X)+\varphi_{i}^{11}(XW))
(\delta_{j}^{11}(X^{-1})+\varphi_{j}^{11}(X^{-1}X_{0}-WY^{-1}))\\\\&&
(\delta_{i}^{11}(X^{-1})+\varphi_{i}^{11}(X^{-1}X_{0}-WY^{-1}))(\delta_{j}^{11}(X)+\varphi_{j}^{11}(XW))],
\end{array}
\end{equation}

\begin{equation}\label{10}\begin{array}{rcl}& &
\varphi_{n}^{12}(X_{0}+X^{-1}X_{0}Y)=\sum\limits_{i+j=n}[(\delta_{i}^{11}(X)+\varphi_{i}^{11}(XW))
(\delta_{j}^{12}(X^{-1})+\varphi_{j}^{12}(X^{-1}X_{0}-WY^{-1})+\tau_{j}^{12}(Y^{-1}))\\\\&&
+(\delta_{i}^{12}(X)+\varphi_{i}^{12}(XW)+\tau_{i}^{12}(Y))(\tau_{j}^{22}(Y^{-1})+\varphi_{j}^{22}(X^{-1}X_{0}-WY^{-1}))
\\\\&&+(\delta_{i}^{11}(X^{-1})+\varphi_{i}^{11}(X^{-1}X_{0}-WY^{-1}))(\delta_{j}^{12}(X)+\varphi_{j}^{12}(XW)+\tau_{j}^{12}(Y))
\\\\&&+(\delta_{i}^{12}(X^{-1})+\varphi_{i}^{12}(X^{-1}X_{0}-WY^{-1})+\tau_{i}^{12}(Y^{-1}))(\tau_{j}^{22}(Y)+\varphi_{j}^{22}(XW)
)],
\end{array}
\end{equation}

\begin{equation}\label{10}\begin{array}{rcl}& &
\varphi_{n}^{22}(X^{-1}X_{0}Y)=\sum\limits_{i+j=n}[(\tau_{i}^{22}(Y)+\varphi_{i}^{22}(XW))
(\tau_{j}^{22}(Y^{-1})+\varphi_{j}^{22}(X^{-1}X_{0}-WY^{-1}))\\\\&&+
(\tau_{i}^{22}(Y^{-1})+\varphi_{i}^{22}(X^{-1}X_{0}-WY^{-1}))(\tau_{j}^{22}(Y)+\varphi_{j}^{22}(XW))].
\end{array}
\end{equation}

Now we take $X=2I_{1}$ and $Y=I_{2}$ in Eq. (44) and Eq. (46), it is
checked that
$$\begin{array}{rcl}& &
\frac{1}{2}\varphi_{n}^{11}(X_{0})=\sum\limits_{i+j=n}[(2\delta_{i}^{11}(I_{1})+2\varphi_{i}^{11}(W))
(\frac{1}{2}\delta_{j}^{11}(I_{1})+\varphi_{j}^{11}(\frac{1}{2}X_{0}-W))\\\\&&
(\frac{1}{2}\delta_{i}^{11}(I_{1})+\varphi_{i}^{11}(\frac{1}{2}X_{0}-W))(2\delta_{j}^{11}(I_{1})+2\varphi_{j}^{11}(W))],
\end{array}
$$\\
$$\begin{array}{rcl}& &
\frac{1}{2}\varphi_{n}^{22}(X_{0})=\sum\limits_{i+j=n}[(\tau_{i}^{22}(I_{2})+2\varphi_{i}^{22}(W))
(\tau_{j}^{22}(I_{2})+\varphi_{j}^{22}(\frac{1}{2}X_{0}-W))\\\\&&+
(\tau_{i}^{22}(I_{2})+\varphi_{i}^{22}(\frac{1}{2}X_{0}-W))(\tau_{j}^{22}(I_{2})+2\varphi_{j}^{22}(W))].
\end{array}
$$\\\\
By the fact that $\delta_{n}^{11}(I_{1})=0 (n\geq 1)$,
$\tau_{n}^{22}(I_{2})=0 (n\geq 1)$ and $\varphi_{n}^{11}(X_{0})=0$,
$\varphi_{n}^{22}(X_{0})=0 $ for any $n\geq 0$, it follows that
$$0=2\varphi_{n}^{11}(W)+4\sum\limits_{i+j=n}\varphi_{i}^{11}(W)\varphi_{j}^{11}(W),$$

$$0=2\varphi_{n}^{22}(W)+4\sum\limits_{i+j=n}\varphi_{i}^{22}(W)\varphi_{j}^{22}(W).$$

When $n=0$, $\varphi_{0}^{11}(W)=\varphi_{0}^{22}(W)=0$, When $n=1$,
$\varphi_{1}^{11}(W)=\varphi_{1}^{22}(W)=0$, So we assume that
$\varphi_{m}^{11}(W)=\varphi_{m}^{22}(W)=0$ for all $1\leq m<n$ and
$W\in \mathcal{M}$. Combining the above equation with the
assumption, we get that $\varphi_{n}^{11}(W)=\varphi_{n}^{22}(W)=0$
for all $1\leq m<n$.

By setting respectively $Y=\frac{1}{2}I_{2}$ and $Y=I_{2}$ in Eq.
(45), the following two equations hold
\begin{equation}\label{10}\begin{array}{rcl}& &
\varphi_{n}^{12}(X_{0}+\frac{1}{2}X^{-1}X_{0})=\sum\limits_{i+j=n}
[\delta_{i}^{11}(X)(\delta_{j}^{12}(X^{-1})+\varphi_{j}^{12}(X^{-1}X_{0}-2W)+2\tau_{j}^{12}(I_{2}))\\\\&&
+\delta_{i}^{11}(X^{-1})(\delta_{j}^{12}(X)+\varphi_{j}^{12}(XW)+\frac{1}{2}\tau_{j}^{12}(I_{2}))]
+2\delta_{n}^{12}(X)\\\\&&
+2\varphi_{n}^{12}(XW)+\tau_{n}^{12}(I_{2})+\frac{1}{2}\delta_{n}^{11}(X^{-1})
+\frac{1}{2}\varphi_{n}^{12}(X^{-1}X_{0}-2W)+\tau_{n}^{12}(I_{2}),
\end{array}
\end{equation}

\begin{equation}\label{10}\begin{array}{rcl}& &
\varphi_{n}^{12}(X_{0}+X^{-1}X_{0})=\sum\limits_{i+j=n}
[\delta_{i}^{11}(X)(\delta_{j}^{12}(X^{-1})+\varphi_{j}^{12}(X^{-1}X_{0}-W)+\tau_{j}^{12}(I_{2}))\\\\&&
+\delta_{i}^{11}(X^{-1})(\delta_{j}^{12}(X)+\varphi_{j}^{12}(XW)+\tau_{j}^{12}(I_{2}))]
+\delta_{n}^{12}(X)\\\\&&
+\varphi_{n}^{12}(XW)+\tau_{n}^{12}(I_{2})+\delta_{n}^{11}(X^{-1})
+\varphi_{n}^{12}(X^{-1}X_{0}-W)+\tau_{n}^{12}(I_{2}).
\end{array}
\end{equation}
Which implies that
\begin{equation}\label{10}\begin{array}{rcl}& &
-\frac{1}{2}\varphi_{n}^{12}(X^{-1}X_{0})=\sum\limits_{i+j=n}
[-\delta_{i}^{11}(X)\varphi_{j}^{12}(W)+\delta_{i}^{11}(X)\tau_{j}^{12}(I_{2})\\\\&&
+\frac{1}{2}\delta_{i}^{11}(X^{-1})\tau_{j}^{12}(I_{2})]
+\delta_{n}^{12}(X)\\\\&&
+\varphi_{n}^{12}(XW)-\frac{1}{2}\delta_{n}^{11}(X^{-1})
-\frac{1}{2}\varphi_{n}^{12}(X^{-1}X_{0}).
\end{array}
\end{equation}

It follows from Eq. (34) and the fact
$\delta_{n}^{22}(X)=\varphi_{n}^{11}(W)=\varphi_{n}^{22}(W)=0$, we
have
\begin{equation}
\delta_{n}^{12}(X)=-\sum\limits_{i+j=n}\delta_{i}^{11}(X)\tau_{j}^{12}(I_{2}).
\end{equation}
Hence combing  Eq. (49) with Eq. (50), we can see that
$$\varphi_{n}^{12}(XW)=\sum\limits_{i+j=n}\delta_{i}^{11}(X)\varphi_{j}^{12}(W)$$
for any invertible $X\in\mathcal{A}$. There exists some $n\in N$
such that $nI_{1}-X$ is invertible for every $X\in \mathcal{A}$, one
can check that
\begin{equation}
\varphi_{n}^{12}(XW)=\sum\limits_{i+j=n}\delta_{i}^{11}(X)\varphi_{j}^{12}(W)
\end{equation}
for any $X\in\mathcal{A}$.

Now we take respectively $X=I_{1}$ and $X=2I_{1}$ in Eq. (45), one
gets
\begin{equation}\label{10}\begin{array}{rcl}& &
\varphi_{n}^{12}(X_{0}+X_{0}Y)=\sum\limits_{i+j=n}[(\delta_{i}^{12}(I_{1})+\varphi_{i}^{12}(W)+\tau_{i}^{12}(Y))\tau_{j}^{22}(Y^{-1})
\\\\&&+(\delta_{i}^{12}(I_{1})+\varphi_{i}^{12}(X_{0}-WY^{-1})+\tau_{i}^{12}(Y^{-1}))\tau_{j}^{22}(Y)]+\delta_{n}^{12}(I_{1})
\\\\&&+\varphi_{n}^{12}(X_{0}-WY^{-1})+\tau_{n}^{12}(Y^{-1})+\delta_{n}^{12}(I_{1})+
\tau_{n}^{12}(Y)+\varphi_{n}^{12}(W),
\end{array}
\end{equation}\\

\begin{equation}\label{10}\begin{array}{rcl}& &
\varphi_{n}^{12}(X_{0}+\frac{1}{2}X_{0}Y)=\sum\limits_{i+j=n}[(2\delta_{i}^{12}(I_{1})+2\varphi_{i}^{12}(W)+\tau_{i}^{12}(Y))\tau_{j}^{22}(Y^{-1})
\\\\&&+(\frac{1}{2}\delta_{i}^{12}(I_{1})+\varphi_{i}^{12}(\frac{1}{2}X_{0}-WY^{-1})+\tau_{i}^{12}(Y^{-1}))\tau_{j}^{22}(Y)]+\delta_{n}^{12}(I_{1})
\\\\&&+2\varphi_{n}^{12}(\frac{1}{2}X_{0}-WY^{-1})+2\tau_{n}^{12}(Y^{-1})+\delta_{n}^{12}(I_{1})+
\frac{1}{2}\tau_{n}^{12}(Y)+\varphi_{n}^{12}(W),
\end{array}
\end{equation}
which implies that
\begin{equation}\label{10}\begin{array}{rcl}& &
\frac{1}{2}\varphi_{n}^{12}(X_{0}Y)=\sum\limits_{i+j=n}[-(\delta_{i}^{12}(I_{1})+\varphi_{i}^{12}(W))\tau_{j}^{22}(Y^{-1})
\\\\&&+\frac{1}{2}(\delta_{i}^{12}(I_{1})+\varphi_{i}^{12}(X_{0}))\tau_{j}^{22}(Y)]
+\varphi_{n}^{12}(WY^{-1})-\tau_{n}^{12}(Y^{-1})+\frac{1}{2}\tau_{n}^{12}(Y).
\end{array}
\end{equation}
Combining the above equation with Eq. (37) and the fact
$\tau_{n}^{11}(Y)=0$, we get
\begin{equation}\label{10}\begin{array}{rcl}& &
\sum\limits_{i+j=n}[-\delta_{i}^{12}(I_{1})\tau_{j}^{22}(Y^{-1})+\frac{1}{2}\delta_{i}^{12}(I_{1})\tau_{j}^{22}(Y)+\frac{1}{2}\varphi_{i}^{12}(X_{0})\tau_{j}^{22}(Y)]+
\frac{1}{2}\tau_{n}^{12}(Y)-\tau_{n}^{12}(Y^{-1})\\\\&&
=\sum\limits_{i+j=n}[-\delta_{i}^{12}(I_{1})\tau_{j}^{22}(Y^{-1})+\frac{1}{2}\delta_{i}^{12}(I_{1})\tau_{j}^{22}(Y)+\frac{1}{2}\varphi_{i}^{12}(X_{0})\tau_{j}^{22}(Y)]
\\\\&&-\sum\limits_{i+j=n}\varphi_{i}^{12}(W)\tau_{j}^{22}(Y^{-1})+ \frac{1}{2}\tau_{n}^{12}(Y)-\tau_{n}^{12}(Y^{-1})+
\varphi_{n}^{12}(WY^{-1}).
\end{array}
\end{equation}
So
\begin{equation}
\varphi_{n}^{12}(WY^{-1})=\sum\limits_{i+j=n}\varphi_{i}^{12}(W)\tau_{j}^{22}(Y^{-1}).
\end{equation}
Replacing $Y$ by $Y^{-1}$ in the above equation, we obtain for any
invertible $Y\in \mathcal{B}$
\begin{equation}
\varphi_{n}^{12}(WY)=\sum\limits_{i+j=n}\varphi_{i}^{12}(W)\tau_{j}^{22}(Y).
\end{equation}
Since there is some integer $n$ such that $nI_{2}-Y$ is invertible
for every $Y\in \mathcal{B}$, it is easy to see that Eq. (57) is
true for every $Y\in \mathcal{B}$ and  $W\in \mathcal{M}$, Summing
up Eq. (54) and Eq. (56), we obtain that
\begin{equation}
\sum\limits_{i+j=n}\delta_{i}^{12}(I_{1})\tau_{j}^{22}(Y^{-1})+\tau_{n}^{12}(Y^{-1})
=\frac{1}{2}[\sum\limits_{i+j=n}\delta_{i}^{12}(I_{1})\tau_{j}^{22}(Y)+\tau_{n}^{12}(Y)].
\end{equation}
Thus
\begin{equation}
\sum\limits_{i+j=n}\delta_{i}^{12}(I_{1})\tau_{j}^{22}(Y)+\tau_{n}^{12}(Y)
=\frac{1}{2}[\sum\limits_{i+j=n}\delta_{i}^{12}(I_{1})\tau_{j}^{22}(Y^{-1})+\tau_{n}^{12}(Y^{-1})]
\end{equation}
by replacing $Y^{-1}$ by $Y$ in the Eq. (58). Combining Eq. (58)
with Eq. (59), we can obtain
$$\frac{1}{2}[\sum\limits_{i+j=n}\delta_{i}^{12}(I_{1})\tau_{j}^{22}(Y^{-1})+\tau_{n}^{12}(Y^{-1})]=2[\sum\limits_{i+j=n}\delta_{i}^{12}(I_{1})\tau_{j}^{22}(Y^{-1})+\tau_{n}^{12}(Y^{-1})].$$
So using the direct computation, we can claim that
\begin{equation}
\tau_{n}^{12}(Y)=-\sum\limits_{i+j=n}
\delta_{i}^{12}(I_{1})\tau_{j}^{22}(Y).
\end{equation}

Now summing up all the above equations and using similar arguments
as that in the proof of Theorem 2.1, it is easily checked that both
$\{\delta_{n}^{11}\}_{n\in N}$ and $\{\tau_{n}^{22}\}_{n\in N}$ are
higher derivations. Therefore it is also an easy computation to see
that $\{D_{n}\}_{n\in N}$ is a higher derivation.
$\Box$\\
\section*{Reference}\small
\begin{description}
  \item[[1]] R.L. An, J.C. Hou, Characterization of derivations on triangular rings: Additive maps derivable at idempotents, 431 (2009) 1070-1080.
  \item[[2]] R.L. An, J.C. Hou, Additivity of Jordan multiplicative maps on Jordan operator algebras. Taiwanese J. Math. 10 (2006) 45--64.
  \item[[3]] W. Jing, On Jordan all-derivable points of
  $\mathcal{B}(\mathcal{H})$, Linear Algebra Appl. 430 (2009) 941-946.
  \item[[4]] Z.K. Xiao, F. Wei, Higher derivations of triangular
  algebras and its generations, Linear Algebra Appl. 432(2010) 2615-2622.
  \item[[5]] J. Zhang, W. Yu, Jordan derivations of triangular
  algebras, Linear Algebra Appl. 419 (2006) 251-255.
  \item[[6]] S. Zhao, J. Zhu, Jordan all-derivable points in the algebra of all upper triangular matrices, Linear Algebra Appl. 433 (2010) 1922-1938.
  \item[[7]] J. Zhu, All-derivable points of operator algebras, Linear Algebra Appl.
  427 (2007) 1-5.

  \end{description}
\end{document}